\documentclass[11pt]{amsart}

\usepackage{amssymb}
\usepackage{caption}
\usepackage{subcaption}
\usepackage{graphicx}
\usepackage{hyperref}
\usepackage{xcolor}
\usepackage{fancyhdr}
\usepackage[normalem]{ulem}
\usepackage[shortlabels]{enumitem}
\usepackage{tikz}

\newcommand{\defn}[1]{{\color{blue} \it {#1}}}

\newcommand{\ZZ}{\mathbb{Z}}

\newcommand{\FF}[1]{\mathbb{F}_{#1}}
\newcommand{\wt}{\mathrm{weight}}
\newcommand{\Min}{M}
\newcommand{\Max}{N}
\newcommand{\Rank}{D}
\newcommand{\Defect}{K}
\newcommand{\graph}{G}
\DeclareRobustCommand{\qbinom}{\genfrac{\lbrack}{\rbrack}{0pt}{}}

\DeclareMathOperator{\supp}{supp}
\DeclareMathOperator{\ZG}{\mathsf{ZG}}
\DeclareMathOperator{\SG}{\mathsf{SG}}
\DeclareMathOperator{\WRG}{\mathsf{WRG}}
\DeclareMathOperator{\WCG}{\mathsf{WCG}}
\DeclareMathOperator{\egraph}{\mathsf{E}}
\DeclareMathOperator{\maxhit}{maxhit}

\usepackage{array}
\newcolumntype{P}[1]{>{\centering\arraybackslash}p{#1}}
\usepackage{multirow}
\usepackage{tabularx}

\numberwithin{equation}{section}
\theoremstyle{definition}
\newtheorem{theorem}{Theorem}[section]
\newtheorem{definition}[theorem]{Definition}
\newtheorem{conjecture}[theorem]{Conjecture}
\newtheorem*{conjecture*}{Conjecture}
\newtheorem{lemma}[theorem]{Lemma}
\newtheorem{proposition}[theorem]{Proposition}
\newtheorem{corollary}[theorem]{Corollary}
\newtheorem*{theorem*}{Theorem}
\theoremstyle{remark}
\newtheorem{remark}[theorem]{Remark}
\newtheorem{example}[theorem]{Example}

\usepackage[margin=1in]{geometry}

\tikzstyle{rednode} = [circle, draw=red!80, fill=red!80, inner sep = 0pt, minimum size = 1.5mm]
\tikzstyle{bluenode} = [circle, draw=blue!80, fill=blue!80, inner sep = 0pt, minimum size = 1.5mm]
\tikzstyle{dot} = [circle, draw=black, fill=black, inner sep = 0pt, minimum size = 0.5mm]

\usepackage[
backend=biber,
style=alphabetic,
sorting=ynt,
maxbibnames=99
]{biblatex}
\title[Positivity properties of $q$-hit numbers]{Positivity properties of $q$-hit numbers\\in the finite general linear group}

\author{Jeffrey Chen}
\address{Massachusetts Institute of Technology}
\email{jeffreychen51@gmail.com}

\author{Jesse Selover}
\address{Department of Mathematics and Statistics\\ University of Massachusetts Amherst\\ Amherst, MA, 01003}
\email{jselover@umass.edu}

\begin{document}
\maketitle
\begin{abstract}
  We consider the problem of counting matrices over a finite field with fixed
  rank and support contained in a fixed set. The count of such matrices gives a
  $q$-analogue of the classical rook and hit numbers, known as the $q$-rook and
  $q$-hit numbers. They are known not to be polynomial in $q$ in general. We use
  inclusion-exclusion on the support of the matrices and the orbit counting
  method of Lewis et al. to show that the residues of these functions in low
  degrees are polynomial. We define a generalization of the classical rook and
  hit numbers which count placements of certain classes of graphs. These give us a
  formula for residues of the $q$-rook and $q$-hit numbers in low degrees. We
  analyze the residues of the $q$-hit number and show that the coefficient of
  $q-1$ in the $q$-hit number is always non-negative.
\end{abstract}
\section{Introduction}

The rook and hit numbers, first studied by Kaplansky and Riordan in
\cite{kaplansky_problem_1946}, count placements of non-attacking rooks on
a subset of the squares of a chess board. These numbers were given $q$-analogues
by Garsia and Remmel in \cite{garsia_q-counting_1986} when the subset is a
Ferrers board. The Garsia-Remmel $q$-rook and $q$-hit numbers are polynomials in
$q$ with non-negative coefficients. Another $q$-rook number was given in
\cite{lewis_matrices_2011}, defined by counting the number of matrices over a
finite field $\mathbb{F}_q$ of a given rank and with zeroes in prescribed
positions (outside the ``board''). This $q$-rook number is sometimes
non-polynomial in $q$, as shown by J.R. Stembridge
\cite{stembridge_counting_1998}. The two definitions of a $q$-rook
number coincide for Ferrers boards by work of Haglund in
\cite{haglund_q-rook_1997}.
The \(q\)-rook numbers of Ferrers boards have been
fruitfully studied from the perspective of coding theory, for example in
\cite{ravagnani_codes_2015}, and the \(q\)-hit numbers of Ferrers boards appear
in the study of the chromatic quasisymmetric function of indifference graphs
\cite{an_chromatic_2021, ls_explicit_2022, nt_down_2022, colmenarejo_chromatic_2023}.

In \cite{lewis_rook_2020}, J.B. Lewis and A.H.
Morales gave a $q$-analogue of the hit numbers for all boards, defined via
inclusion-exclusion in terms of the Lewis et al. $q$-rook numbers. They showed
that when the board is a Ferrers board this definition coincides with Garsia and
Remmel's $q$-hit numbers. However, as with the Lewis et al. $q$-rook numbers, these
$q$-hit numbers are not always polynomial. They are conjectured to have many
strong non-negativity properties, including being non-negative for fixed $q$ and
any board.

Consider a \defn{board} \(B\), a subset of the grid \([m] \times [n]\). We can
think of \(B\) as an \(m \times n\) matrix where \(B_{ij}\) is \(1\) if the cell
\((i,j)\) is in \(B\) and \(0\) otherwise. Define the \defn{rook number}
\(r_i(B)\) as
the number of ways to place $i$ rooks on the cells of $B$ such that no two
attack each other. Also define the \defn{hit number} $h_i(B)$ as the number of ways to
place $m$ non-attacking rooks on $[m]\times[n]$ with exactly $i$ rooks in $B$
\cite{kaplansky_problem_1946}. They are related by the following equation:

\begin{equation}
\label{eqn:hit-rook-relation}
\sum_{i=0}^{n} h_i(B) t^i = \sum_{i=0}^{n} r_i(B) \frac{(n-i)!}{(n-m)!}(t-1)^{i}.
\end{equation}

Define the \defn{support} of a matrix $A$ to be the board \(\supp(A)\) with
\((i,j) \in \supp(A)\) if $A_{i,j}$ is nonzero. For a board
$B\subseteq[m]\times[n]$, define the \defn{matrix count} $\mathfrak{m}_i(B,q)$ as the number of $m$ by
$n$ matrices in $\FF{q}$ with rank \(i\) and support contained in $B$.
Lewis et al.\cite[Prop. 5.1]{lewis_matrices_2011} showed that

\[\mathfrak{m}_i (B,q)\equiv r_i (B) (q-1)^{i} \mod (q-1)^{i+1}.\]

This means that for a fixed $q$, $\mathfrak{m}_i (B,q)$ is
divisible by $(q-1)^{i}$. Define the \defn{\(q\)-rook number}
$M_i (B,q) := \mathfrak{m}_i (B,q) / (q-1)^i$.
Thus $M_i(B,q)$ is an integer for
fixed $q$, and since
$M_i(B,q) \equiv r_i(B) \pmod{q-1}$, $M_i(B,q)$ is a $q$-analogue of
the rook numbers $r_i(B)$. We know that $\mathfrak{m}_i (B,q)$ is not
necessarily a polynomial in $q$ \cite{stembridge_counting_1998}. However, for certain classes of boards, this
count has been proven to be polynomial \cite{stembridge_counting_1998}[Theorem
8.2]
, and for these classes of board \(M_i(B,q)\) necessarily is polynomial as well.

Lewis and Morales defined a $q$-analogue of the hit numbers in \cite[Eq.
(1.3)]{lewis_rook_2020}, by the generating function equality

\begin{equation}
\sum_{i=0}^{n}H_i(B,q)t^i = q^{\binom{m}{2}} \sum_{i=0}^n M_i(B,q) \cdot \frac{[n-i]!_q}{[n-m]!_q} (t-1)(tq^{-1} - 1)\ldots (tq^{-(i-1)} - 1),
\end{equation}
a \(q\)-analogue of Equation~\eqref{eqn:hit-rook-relation}.

This paper was motivated by conjectures in the paper \cite{lewis_rook_2020} by
Lewis and Morales. Firstly, that \(H_d(B,q)\) is non-negative for any board
\(B\) and prime power \(q\) (\cite[Conjecture 6.3]{lewis_rook_2020}), and
secondly, that for a certain class of boards \(B\) where \(H_d(B,q)\) is
polynomial, the polynomial \(H_d(B,x+1)\) has non-negative coefficients in \(x\)
(\cite[Conjecture 6.7]{lewis_rook_2020}).

We find a positivity result for the coefficients of \(H_d(B, x+1)\) which applies unconditionally to all boards, and so
we must say what we mean by the ``coefficient'' when this function of \(x\) is non-polynomial.

Let \(S\) be an unbounded subset of
\(\mathbb{Z}\), and let \(f\) be a function \(S \to \mathbb{Q}\). When there is
a polynomial \(P \in \mathbb{Z}[x]\) of degree at most \(k-1\) such that
\[
P(x) \equiv f(x) \pmod{x^k}
\]
for all \(x\) in \(S\), we say that \(f\) is \defn{polynomial modulo} \(x^k\),
and we say the \defn{coefficient} of \(x^i\) in \(f\) is the coefficient of
\(x^i\) in \(P\) for \(i < k\). Note that there is at most one polynomial \(P\)
satisfying these conditions for any fixed \(k\), so this coefficient is
well-defined. If \(R\) is an unbounded subset of \(S\) and \(f|_R\) is
polynomial modulo \(x^k\), we say that \(f\) is polynomial modulo \(x^k\)
\emph{on \(R\)}, and we may talk about \(f\)'s coefficients on \(R\).

Our main result is
\begin{theorem*}[Theorem \ref{thm:q-coeff1-pos}]
  For a board $B \subseteq [n]\times[n]$, \(H_d(B,x+1)\) is polynomial modulo
  \(x^2\), and its coefficient of \(x\) is non-negative.
\end{theorem*}

In Section~\ref{sec:q-hit-poly} we extend the polynomiality result by showing that \(\mathfrak{m}_i(B,x+1)\),
\(M_i(B, x+1)\), and \(H_i(B, x+1)\) are polynomial modulo \(x^6\) for any
board. Neither the positivity, nor the polynomiality, extend to the functions
\(\mathfrak{m}_i(B,q), M_i(B,q)\), or \(H_i(B,q)\) modulo powers of \(q\) for
all boards.

In fact, though, we know of no counterexample to the stronger conjecture:

\begin{conjecture*}[Conjecture~\ref{conj:q-poly-all-coeffs}]
  Let \(B\) be a board, and let \(S\) be an
  unbounded subset of \(\mathbb{Z}\) such that \(H_d(B, x+1)\) is polynomial
  modulo \(x^k\) on \(S\) for some \(k\). Then the coefficients of \(x^i\), for
  \(i < k\), in \(H_d(B, x+1)\) on \(S\) are non-negative.
\end{conjecture*}

If \(B\) is a board such that \(H_d(B, x+1) = P(x)\) is polynomial, then
\(H_d(B, x+1)\) is polynomial modulo \(x^k\) for all \(k\), and the coefficients
of \(H_d(B, x+1)\) modulo \(x^k\) agree with the first \(k\) coefficients of
\(P\). Thus Conjecture~\ref{conj:q-poly-all-coeffs} would imply Conjecture~6.3
and Conjecture~6.7 in \cite{lewis_rook_2020}.

We note that \(H_d(B, x+1)\) is suspected to be
``polynomial on residue classes,'' i.e. polynomial on the set
\(\{x \mid x \equiv a \pmod{k}\}\) for some \(k\) and each \(a\) (as in Example~\ref{exa:fano}).
For a board \(B\) where \(H_d(B, x+1)\) is congruent to \(f_a(x+1)\) on the residue
class \(\{x \mid x \equiv a \pmod{k}\}\), our
Conjecture~\ref{conj:q-poly-all-coeffs} says that \(f_a(x+1)\) is a polynomial
with non-negative integer coefficients.

The \(k=1\) case of Conjecture~\ref{conj:q-poly-all-coeffs} is Proposition~3.3
in \cite{lewis_rook_2020}. Our Theorem~\ref{thm:q-coeff1-pos} is the \(k=2\)
case. Depending on your viewpoint, it may either be partial evidence for Lewis
and Morales's conjecture \cite[Conjecture 6.7]{lewis_rook_2020}, or constrain
the search for counterexamples.

We also give a possible framework for proving higher-$k$ cases, although
we do not expect our methods to extend past $k=6$.

Our proof technique is explained in Section~\ref{sec:orbit-q-rook}, where we
define generalized rook and hit numbers (Definition~\ref{def:gen-rook}) that we
hope are of independent interest.
We give signed formulas in terms of these generalized rook and hit numbers for
the coefficients of \(M_i(B, x+1)\) modulo \(x^2\) (Theorem~\ref{thm:q-rookres})
and \(H_i(B, x+1)\) modulo \(x^2\) (Theorem~\ref{thm:q-hitres}). In Section~
\ref{sec:coeff-hit-positive} we prove the inequalities that allow us to show the
coefficient of \(x\) in \(H_i(B, x+1)\) is non-negative.

\noindent
\textbf{Outline.}
In Section~\ref{sec:background}, we give definitions of the
$q$-rook and $q$-hit numbers and recall important results from the literature.
In Section~\ref{sec:orbits}, we give a formula for $M_d(B,x+1) \pmod{x^{2}}$
in terms of generalized rook numbers of $B$. We also study the refinement of the
set of matrices with rank $d$ and support in $B$ to show that
$M_{d}(B,x+1)$  must be polynomial modulo \(x^6\) for any board $B$. In
Section~\ref{sec:higher-hit}, we transform our formula for
$M_d(B,x+1) \pmod{x^2}$ into a formula for \(H_d(B,x+1) \pmod{x^2}\).
We show that the coefficients of $H_d(B,x+1)\pmod{x^2}$ are non-negative.

\section{Background information}
\label{sec:background}

In this section, we give the definitions and background information about the
$q$-rook and $q$-hit numbers, and then we review important past results.

First, consider the $q$-analogues of the natural numbers 
\[[i]_q = (q^{i} - 1)/(q-1) = 1 + q + \cdots + q^{i-1}.\]
We define the \defn{\(q\)-factorial} and \defn{\(q\)-binomial}
\[[n]!_q = [n]_q[n-1]_q\ldots [1]_q\]
and
\[\qbinom{n}{m}_q = \frac{[n]!_q}{[m]!_q [n-m]!_q}.\]
Define the \defn{$q$-Pochhammer symbol} as
\[(t;q)_k = \prod_{i=0}^{k-1} (1 - tq^{i}) = (1-t)(1-tq)\ldots (1-tq^{k-1}),\]
and define the \defn{$q$-hit numbers} for a board $B\subseteq [m]\times[n]$, where
$m \leq n$, with the equation \cite{lewis_rook_2020}:

\begin{equation}
\label{eqn:q-hitdef}
\sum_{i=0}^{n} H_i(B,q) t^i = q^{\binom{m}{2}} \sum_{i=0}^{m} M_i(B,q) \frac{[n-i]!_q}{[n-m]!_q} (-1)^i (t;q^{-1})_i.
\end{equation}

The $q$-hit numbers are a $q$-analogue in the sense shown by Lewis and Morales
in \cite[Prop. 3.3]{lewis_rook_2020}, which states that
$H_d(B,q) \equiv h_d(B) \pmod{q-1}$. When \(q\) is not \(1\), there is no known
combinatorial interpretation of $H_d(B,q)$ for general boards \(B\). However,
the numbers are conjectured to be positive by Lewis and Morales in
\cite[Conjecture 6.3]{lewis_rook_2020}.

Lewis and Morales also showed the following.

\begin{proposition}[{\cite[Prop. 3.5]{lewis_rook_2020}}]
\label{prop:q-hit-direct}
For any board $B$ inside \([m]\times [n]\), we can compute the individual $q$-hit and
$q$-rook numbers in the following way:
\[H_k (B,q) = q^{\binom{k+1}{2} + \binom{m}{2}} \sum_{i=k}^{m} M_i (B,q) \cdot \frac{[n-i]!_q}{[n-m]!_q} \qbinom{i}{k}_q (-1)^{i+k} q^{-ik}\]
and 
\[M_{k} (B,q) = q^{\binom{k}{2} - \binom{m}{2}} \frac{[n-m]!_q}{[n-k]!_q} \sum_{i=k}^{m} H_{i}(B,q) \qbinom{i}{k}_q.\]
\end{proposition}
\label{prop:LM17-3.5}

By this proposition, if $M_i(B,q)$ is an integer polynomial in $q$ for all \(i\)
such that
$k\leq i \leq m$, then $H_k(B,q) \in \ZZ[q]$, and vice versa. In fact, this
implies that all \(M_i(B,q)\) are polynomials in \(q\) with integer coefficients
if and only if the same is true for all \(H_i(B,q)\).

We also have $H_m(B,q) = M_m(B,q)$ \cite[Prop. 3.5]{lewis_rook_2020}.

\begin{example}
\label{exa:fano}
It is important to note that $H_d(B,q)$ is not a polynomial in $q$ for all
choices of $B$. A counterexample to this is the Fano board
$F \subseteq [7]\times[7]$:

\begin{center}
{\small
$\begin{bmatrix}
*&*&0&0&0&*&0\\
0&*&*&0&0&0&*\\
*&0&*&*&0&0&0\\
0&*&0&*&*&0&0\\
0&0&*&0&*&*&0\\
0&0&0&*&0&*&*\\
*&0&0&0&*&0&*
\end{bmatrix}.$
}
\end{center}
In \cite{stembridge_counting_1998}, Stembridge found that \(M_7(F,q)\) is the
following non-polynomial function of \(q\):
\begin{equation}
\label{eqn:fano-hit-eq}
\begin{split}
M_7(F,q) =& (x+1)^3 (x^{11}+ 17x^{10} + 135x^9 + 650x^8 + 2043x^7 + 4236x^6 \\
+& 5845x^5 + 5386x^4 + 3260x^3 + 1236x^2 + 264x + 24 - Z_2x^6),
\end{split}
\end{equation}
where $x = q-1$ and $Z_2 = 1$ for odd $q$ and $0$ for even $q$.
\end{example}
This shows that the $q$-rook number $M_d(B,q)$ is not always a polynomial in
$q$, and consequently that \(H_d(B,q)\) is not either. See \cite{stanley_spanning_1998} for further discussion on
non-polynomiality of this and related counting problems over $\FF{q}$.

\begin{remark}
\label{rem:fano_minimal}

For all $m$ by $n$ boards $B$ with $m,n\leq 6$, $M_d(B,q)$ is a polynomial in \(q\).
\end{remark}

\section{Orbits of matrices}
\label{sec:orbits}

In this section, we study a group action on matrices with a fixed rank and
support. Following \cite{lewis_matrices_2011}, we find the size of the orbit of
a given matrix with rank $d$. Many orbits have size divisible by
\((q-1)^{d+2}\), and we are able to enumerate the orbits which do not, which
gives a formula for $M_d(B,q) \pmod{(q-1)^2}$. We use the same technique to show
polynomiality of $M_d(B,q)$ modulo $(q-1)^{6}$.

First, we define the graph of a board.

\begin{definition}
  We define a \defn{bi-colored graph} $G$ as a graph together with a partition of its
  vertices into two sets $V_1$ and $V_2$, such that every edge is between $V_1$
  and $V_2$. We call elements of \(V_1\) row vertices, and elements of \(V_2\)
  column vertices. A homomorphism of bi-colored graphs is required to take row
  vertices to row vertices, and column vertices to column vertices. In this
  paper, a graph means a bi-colored graph unless otherwise mentioned.
\end{definition}

When \(B\) is a board (a subset of \([m]\times[n]\)), we define \defn{\graph(B)}
to be the graph with row vertices \([m]\) and column vertices \([n]\), and edge
set \(B\).

\subsection{Counting matrices by support}

We show a relation between the maximal rook placement of a board and the maximal
rank of a matrix with support in that board.

The following is a standard application of Hall's marriage theorem, which we
include for completeness.
\begin{proposition}
\label{prop:max-rank}
Consider a board $B$, and let $k$ be the maximum non-attacking rooks that can be
placed in $B$. For any matrix $M$ over $\FF{q}$ with support in $B$, the rank of
$M$ is at most $k$.
\end{proposition}
\begin{proof}
  Let the rank of an $m$ by $n$ matrix $M$ with support in $B$ be $k$. We show
  that we can place at least $k$ non-attacking rooks on $B$. Let
  $v_1, v_2, \ldots, v_k$ be $k$ rows of $M$ that are linearly independent,
  which exist because the rank is $k$. Now, let $G$ be the bi-colored graph
  formed with $a_1, a_2, \ldots, a_k$ row vertices and $b_1, b_2, \ldots, b_n$
  column vertices, with an edge between $a_i$ and $b_j$ if the $j$th element in
  $v_i$ is nonzero.

  A matching in $G$ with $k$ edges corresponds to $k$ cells in \(B\), no two of
  which are in the same row nor column. By Hall's marriage theorem
  \cite{hall_subsets_1935}, if, for every set $S$ of $a_i$s, there are at least
  $s:=\#S$ nodes incident to some vertex in $S$, then a maximal matching exists.

  For the sake of contradiction, assume there are fewer than $s$ nodes incident
  to at least one vertex in $S$. Without loss of generality, let these nodes be
  $a_1, a_2, \ldots, a_s$ in $S$, and the nodes incident to these be
  $b_1, b_2, \ldots, b_i$ for $i < s$. This corresponds to the first $s$ rows
  only having entries in the first $i$ columns for $i < s$. Therefore the rank
  of $v_1, v_2, \ldots, v_{s}$ is at most $i$, so they are not linearly
  independent, a contradiction.

  We conclude that $M$ has rank at most $k$, where $k$ is the maximum
  non-attacking rooks.
\end{proof}

\begin{definition}
  Define \defn{$S_{d}(B,q)$} as the set of $m$ by $n$ matrices $A$ over $\FF{q}$
  such that the rank of $A$ is $d$ and the support of $A$ is \emph{exactly} $B$.
\end{definition}

Let \defn{$T_q(m,n,B,d)$} be the set of $m$ by $n$ matrices with support
contained in $B$ and rank $d$. We have the following relations:

\[T_q(m,n,B,d) = \bigcup_{C\subseteq B} S_d(C,q)\]

and thus

\[\mathfrak{m}_d(B,q) = \sum_{C\subseteq B} \# S_d(C,q).\]

By M\"obius inversion, we get

\[\# S_d(B,q) = \sum_{C\subseteq B}(-1)^{(\#B - \#C)} \mathfrak{m}_d(C,q).\]

Define \defn{$\maxhit(B)$} as the maximum number of non-attacking rooks that can
be placed on $B$. If $\maxhit(C) < d$, then $\mathfrak{m}_d(C,q) = 0$ by
Proposition~\ref{prop:max-rank}. Thus our two equations become:

\[T_q(B,d) = \bigcup_{C\subseteq B, \maxhit(C) \geq d} S_d(C,q)\]
and

\begin{equation}
\label{eqn:inverted-Sq-formula}
\#S_d(B,q) = \sum_{C\subseteq B, \maxhit(C) \geq d} (-1)^{(\#B-\#C)} \mathfrak{m}_d(C,q).
\end{equation}

\begin{remark}
  By Remark~\ref{rem:fano_minimal} and equation~\eqref{eqn:inverted-Sq-formula},
  the Fano board is also the minimal board such that $\# S_{d}$ is
  non-polynomial.
\end{remark}

\begin{proposition}
\label{prop:orbitSize}
For fixed $q$, the number $\#S_d(B,q)$ is divisible by
$(q-1)^{m + n - C(G(B))}$, where \defn{$C(G(B))$} is the number of connected
components of $G(B)$.
\end{proposition}

\begin{proof}
  We mimic the proof of \cite[Prop. 5.1]{lewis_matrices_2011}. Let
  $A\in S_d(B,q)$ be a matrix. Let $(\FF{q}^{\times})^{l}$ be the set of
  diagonal $l\times l$ matrices with each diagonal entry nonzero. Now, consider
  the group action $(\FF{q}^{\times})^{m} \times (\FF{q}^{\times})^{n}$ on
  $S_{d}$ defined by $(X,Y) \cdot A = XAY^{-1}$. The support of $XAY^{-1}$ is
  still exactly $B$ because $X$ and $Y$ are diagonal matrices. Define
  $x_1, x_2, \ldots, x_m$ and $y_1, y_2, \ldots, y_n$ as the diagonal entries of
  $X$ and $Y$, in that order. We show that $(X,Y)$ stabilizes $A$ if, for each
  nonzero $a_{i,j}$, we have $x_i = y_j$. This is because, if $XAY^{-1} = A$,
  then $XA = AY.$ Then $(XA)_{i,j}$ (the element on the ith row and jth column
  of $XA$) is $x_i a_{i,j}$, and similarly $(AY)_{i,j}$ is $a_{i,j} y_j$. Thus
  if $a_{i,j} \neq 0$, then $x_i = y_j$.

  This means there are $(q-1)^{C(G(B))}$ choices for $X$ and $Y$, because for
  each connected component in $G(B)$, we have $q-1$ ways to choose those
  elements in $X,Y$ over $\FF{q}$. By the orbit-stabilizer theorem, since there
  are $(q-1)^{m+n}$ ways to choose $X$ and $Y$, the size of the orbit of $A$ is
  $(q-1)^{m+n - C(G(B))}$. Finally, this implies that $\# S_{d}(B,q)$ is
  divisible by $(q-1)^{m + n - C(G(B))}$, since we can partition $S_{d}(B,q)$
  into orbits of size $(q-1)^{m+n-C(G(B))}$.
\end{proof}

\subsection{Classifying orbits of size \texorpdfstring{\((q-1)^{d+1}\)}{(q-1)\^(d+1)} in \texorpdfstring{\(\mathfrak{m}_d(B,q)\)}{m\_d(B,q)}}
\label{sec:orbit-q-rook}

Let \(\mathcal{O}\) be the set of orbits of the action of
$(\FF{q}^{\times})^m \times (\FF{q}^{\times})^n$ on \(T_q(m,n,B,d)\),
defined by $(X,Y) \cdot A = XAY^{-1}$. Since the action preserves the support,
there is a well-defined board \(\supp(O)\) for each \(O \in \mathcal{O}\). Therefore,
\[
\mathfrak{m}_d(B,q) = \sum_G \sum_{\substack{B' \subseteq B \\ G(B') \cong G}} \sum_{\substack{O \in \mathcal{O} \\ \supp(O) = B'}} \#O,
\]
where the outer sum is over all isomorphism classes of bi-colored graphs \(G\), but can be restricted to a finite sum over just those isomorphism classes which appear as \(G(B')\) for some \(B' \subseteq B\).

By the orbit-stabilizer theorem,
\[
\#O = (q-1)^{m+n-C(G)},
\]
so
\[
\mathfrak{m}_d(B,q) = \sum_G \sum_{\substack{B' \subseteq B \\ G(B') \cong G}} \sum_{\substack{O \in \mathcal{O} \\ \supp(O) = B'}} (q-1)^{m+n-C(G)},
\]
and grouping terms with the same power of \((q-1)\) we get
\begin{equation}
\label{eqn:m-sum-orbits}
\mathfrak{m}_d(B,q) = \sum_{i=d}^{m+n-1}
(q-1)^i \sum_{\substack{G \\ C(G) = m+n-i}} \sum_{\substack{B' \subseteq B \\ G(B') \cong G}} \sum_{\substack{O \in \mathcal{O} \\ \supp(O) = B'}} 1.
\end{equation}

Furthermore, if \(B_1\) and \(B_2\) are boards such that
\(G(B_1) \cong G(B_2)\), \(S_d(B_1, q)\) has the same number of orbits as
\(S_d(B_2, q)\).

We can rewrite Equation~\eqref{eqn:m-sum-orbits} to remove the inner sums.
Define \defn{\(\mathcal{B}(B, G)\)} to be the number of boards \(B' \subseteq B\) such
that \(G(B') \cong G\),
and \defn{\(\mathcal{O}_d(G,q)\)} to be the number of orbits with support equal
to any one fixed board \(B'\) with \(G(B') \cong G\). Of course,
\(\mathcal{O}_d(G,q)\) is equal to $\#S_d(B',q)/(q-1)^{m+n-C(G)}$ for any such
board \(B'\). If \(G\) has \(x\) row vertices and \(y\) column vertices (fewer than \(m\) and \(n\) respectively), it's impossible for \(G(B')\) to be isomorphic to \(G\), so extend the definitions by setting \(\mathcal{B}(B,G) := \mathcal{B}(B,G \bigsqcup R)\) and \(\mathcal{O}_d(G,q) := \mathcal{O}_d(G \bigsqcup R, q)\) where \(R\) is the completely disconnected graph with \(m-x\) row vertices and \(n-y\) column vertices.

With these definitions, Equation~\eqref{eqn:m-sum-orbits} becomes

\begin{equation}
\label{eqn:m-sum}
    \mathfrak{m}_d(B,q) = \sum_{i=d}^{m+n-1} (q-1)^i
\sum_{\substack{G \\ C(G) = m+n-i}} \mathcal{B}(B, G) \cdot \mathcal{O}_d(G,q),
\end{equation}

One benefit of this expansion is that it is multiplicative:
If \(G_1, G_2\) are bi-colored graphs on different
vertex sets we have
\begin{align*}
  & \mathcal{O}_d(G_1 \sqcup G_2, q)\\
  =& \sum_i \mathcal{O}_i(G_1,q)\mathcal{O}_{d-i}(G_2,q).
\end{align*}
Therefore,
\begin{equation*}
  \mathfrak{m}_d(B,q) = \sum_{i=d}^{m+n-1} (q-1)^i
  \sum_{\substack{G \\ C(G) = m+n-i}} \mathcal{B}(B,G) \sum_{\Rank_1 + \cdots + \Rank_{m+n-i} = d} \prod_{k=1}^{m+n-i} \mathcal{O}_{\Rank_k}(G_k,q),
\end{equation*}
where \(G_1, \ldots, G_{m+n-i}\) are the connected components of \(G\). If \(G_k\) has \(x_k\) row vertices and \(y_k\) column vertices, we can rewrite the equation as
\begin{equation}
  \label{eqn:m-sum-connected}
  \mathfrak{m}_d(B,q) = (q-1)^d \sum_{i=1}^{m+n-d}
  \sum_{\substack{G \\ C(G) = i}} \mathcal{B}(B,G) \sum_{\Rank_1 + \cdots + \Rank_{i} = d} \prod_{k=1}^{i} \mathcal{O}_{\Rank_k}(G_k,q)(q-1)^{x_k+y_k-\Rank_k-1},
\end{equation}

From Equation~\ref{eqn:m-sum-connected}, the single edge graph \defn{\(\egraph\)},
with one row and one column vertex, takes on special importance for low
coefficients of \((q-1)\). The contributions to the \((q-1)^d\) coefficient come from indices where \(\mathcal{O}_{\Rank_k}(G_k,q)\) is nonzero while \(\Rank_k = x_k+y_k-1\) for all \(k\). For each \(k\), this condition can only be satisfied if \(G_k\) is the edge graph \(\egraph\) -- any other graph will have \(\mathcal{O}_{x_k+y_k-1}(G_k,q) = 0\).
In fact, \(\mathcal{O}_1(\egraph,q) = 1\), so the lowest
coefficient of Equation~\ref{eqn:m-sum-connected} is \(\mathcal{B}(B,\egraph^d) = r_d(B)\), the classical rook number
(\cite[Prop 5.1]{lewis_matrices_2011}).

Motivated by this case, we view \(\mathcal{B}(B,G)\)
as a generalized rook number of the board \(B\). Due to the importance of
the single-edge graph for low degree coefficients, we give notation
in Definition~\ref{def:gen-rook} which more closely matches the usual rook numbers.

Modulo any fixed power of \((q-1)\), only finitely many isomorphism classes of
connected bi-colored graphs give nonzero contributions to the inner sum in
Equation~\ref{eqn:m-sum-connected}. To understand \(\mathfrak{m}_d(B,q)\) modulo \((q-1)^{d+2}\), it suffices
to understand the connected bi-colored graphs \(G\) and ranks \(\Rank\) such that
\(x+y-1 = \Rank + 1\) (and \(\mathcal{O}_\Rank(G,q)\) is nonzero), which we enumerate in
Proposition~\ref{prop:graph-classes}.
Then for each graph, and rank, we
calculate \(\mathcal{O}_\Rank(G,q)\) in Proposition~\ref{prop:og}.
The results are gathered into an expansion for the \(q\)-rook number
\(M_d(B,q)\) (modulo \((q-1)^2\)) in Theorem~\ref{thm:q-rookres}.

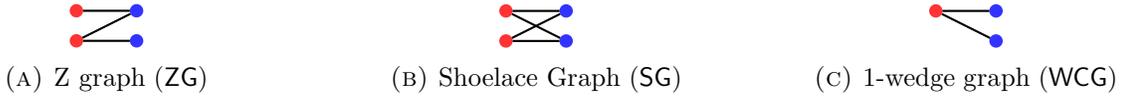
\begin{figure}[htbp]
\centering

\begin{subfigure}[b]{0.3\textwidth}
\centering
\begin{tikzpicture}[x = 0.5cm, y = 0.5cm, thick, rednode/.style = {circle, draw=red!80, fill=red!80, inner sep = 0pt, minimum size = 1.5mm}, bluenode/.style = {circle, draw=blue!80, fill=blue!80, inner sep = 0pt, minimum size = 1.5mm}, dot/.style = {circle, draw=black, fill=black, inner sep = 0pt, minimum size = 0.5mm}, scale=0.8]

\node[rednode] (R2) at (0,0) {};
\node[rednode] (R3) at (0,-1) {};
\node[bluenode] (B2) at (2,0) {};
\node[bluenode] (B3) at (2,-1) {};
\path[-] (R3) edge (B2);
\path[-] (R2) edge (B2);
\path[-] (R3) edge (B3);
\end{tikzpicture}
\caption{Z graph ($\ZG$)}
\label{graph:z}
\end{subfigure}
\hfill
\begin{subfigure}[b]{0.3\textwidth}
\centering
\begin{tikzpicture}[x = 0.5cm, y = 0.5cm, thick, rednode/.style = {circle, draw=red!80, fill=red!80, inner sep = 0pt, minimum size = 1.5mm}, bluenode/.style = {circle, draw=blue!80, fill=blue!80, inner sep = 0pt, minimum size = 1.5mm}, dot/.style = {circle, draw=black, fill=black, inner sep = 0pt, minimum size = 0.5mm}, scale=0.8]

\node[rednode] (R3) at (0,-1) {};
\node[rednode] (R4) at (0,-2) {};
\node[bluenode] (B3) at (2,-1) {};
\node[bluenode] (B4) at (2,-2) {};
\path[-] (R3) edge (B3);
\path[-] (R3) edge (B4);
\path[-] (R4) edge (B3);
\path[-] (R4) edge (B4);
\end{tikzpicture}
\caption{Shoelace Graph ($\SG$)}
\label{graph:shoelace}
\end{subfigure}
\hfill
\begin{subfigure}[b]{0.3\textwidth}
\centering
\begin{tikzpicture}[x = 0.5 cm, y = 0.5cm, thick, rednode/.style = {circle, draw=red!80, fill=red!80, inner sep = 0pt, minimum size = 1.5mm}, bluenode/.style = {circle, draw=blue!80, fill=blue!80, inner sep = 0pt, minimum size = 1.5mm}, dot/.style = {circle, draw=black, fill=black, inner sep = 0pt, minimum size = 0.5mm}, scale=0.8]

\node[rednode] (R3) at (0,-1) {};
\node[bluenode] (B3) at (2,-1) {};
\node[bluenode] (B4) at (2,-2) {};
\path[-] (R3) edge (B3);
\path[-] (R3) edge (B4);
\end{tikzpicture}

\caption{1-wedge graph ($\WCG$)}
\label{graph:1-wedge}
\end{subfigure}
\hfill

\caption{Three bi-colored graphs for $\ZG, \SG, \WRG,$ and $\WCG$}
\label{fig:three-graphs}
\end{figure}

\begin{definition}
\label{def:graph-classes}
We define four bi-colored graphs as follows, illustrated by Figure~\ref{fig:three-graphs}:
\begin{enumerate}
    \item Define \defn{$\ZG$} as the bi-colored graph with $2$ rows and $2$ columns, which forms a path of length $3$ (Figure~\ref{graph:z}).
    \item Define \defn{$\SG$} as the bi-colored graph with $2$ rows and $2$ columns, which forms a $K_{2,2}$ (Figure~\ref{graph:shoelace}).
    \item Define \defn{$\WCG$} as the bi-colored graph with $1$ row and $2$ columns, which forms a row connected to two columns (Figure~\ref{graph:1-wedge}).
    \item Define \defn{$\WRG$} as the bi-colored graph with $2$ rows and $1$ column, which forms a column connected to two rows (Figure~\ref{graph:1-wedge} with red and blue flipped).
    
\end{enumerate}
\end{definition}

\begin{proposition}
\label{prop:graph-classes}
If \(G\) is a connected bi-colored graph with \(x\) row vertices and \(y\) column
vertices, \(\Rank\) is a rank such that \(\mathcal{O}_{\Rank}(G,q)\) is nonzero,
and \(x+y-1 = \Rank+1\), then \((G,\Rank)\) is one of:
\begin{itemize}
\item \((\ZG,2)\)
\item \((\SG,2)\)
\item \((\WCG,1)\)
\item \((\WRG,1)\).
\end{itemize}
\end{proposition}
\begin{proof}
Since \(\Rank \leq \min(x,y)\), and \(x+y = \Rank+2\), we must have
\(\max(x,y) \leq 2\). The proof follows by checking the finite list of bi-colored
graph isomorphism classes with at most two of each color of vertex.
\end{proof}

\begin{proposition}
\label{prop:og}
    The number of orbits of the appropriate rank for each of these graphs is
    \begin{enumerate}
        \item[(a)] \(\mathcal{O}_2(\ZG,q) = 1\).
        \item[(b)] \(\mathcal{O}_2(\SG,q) = q-2\).
        \item[(c)] \(\mathcal{O}_1(\WCG,q) = \mathcal{O}_1(\WRG, q) = 1\).
    \end{enumerate}
\end{proposition}

\begin{figure}[htbp]
\label{fig:three-boards}
\centering

\begin{subfigure}[b]{0.3\textwidth}
\centering
\begin{equation*}
  \left[\begin{matrix}
* & * \\ 0 & *
        \end{matrix} \right]
\end{equation*}
\caption{Board for \(\ZG\)}
\label{fig:board-zg}
\end{subfigure}
\begin{subfigure}[b]{0.3\textwidth}
\centering
\begin{equation*}
  \left[\begin{matrix}
* & * \\ * & *
        \end{matrix} \right]
\end{equation*}
\caption{Board for \(\SG\)}
\label{fig:board-sg}
\end{subfigure}
\begin{subfigure}[b]{0.3\textwidth}
\centering
\begin{equation*}
  \left[\begin{matrix}
* & * \\ 0 & 0
        \end{matrix} \right]
\end{equation*}
\caption{Board for \(\WCG\)}
\label{fig:board-wcg}
\end{subfigure}
\caption{Boards used in the proof of Proposition~\ref{prop:orbitSize}}
\end{figure}

\begin{proof}
  In each case we can calculate \(\#S_\Rank(B,q)\) for some board \(B\) with the
  correct (bi-colored) graph and divide by the orbit size \((q-1)^{x+y-1}\).

  Part (a):
  Let \(B\) be the board from Figure~\ref{fig:board-zg} which has graph
  \(G(B) = \ZG\). Any \(2 \times 2\) matrix with this support has rank \(2\) by
  counting pivots, so there are \((q-1)^3\) matrices in \(S_2(B,q)\), and
  \(\mathcal{O}_2(\ZG,q) = 1\).

  Part (b): Let $B$
  be the board from Figure~\ref{fig:board-sg} which has graph
  \(G(B) = \SG\). There are \((q-1)^4\) matrices with that support over
  \(\FF{q}\), \((q-1)^3\) of which have rank \(1\), so there are
  \((q-1)^3(q-2)\) matrices in \(S_2(B,q)\), and \(\mathcal{O}_2(\SG,q) = q-2\).

  Part (c):
  Let \(B\) be the
  board from Figure~\ref{fig:board-wcg} which has graph
  \(G(B) = \WCG\). Any graph with this support has rank \(1\) by counting
  pivots, so there are \((q-1)^2\) matrices in \(S_1(B,q)\) and
  \(\mathcal{O}_1(\WCG,q) = 1\). The case of \(\WRG\) is symmetric.
\end{proof}

For convenience, we define generalized rook numbers:
\begin{definition}
\label{def:gen-rook}
For a board $B\subseteq[m]\times[n]$ and a bi-colored graph $F$ with $x$ row and
$y$ column vertices, we define \defn{$r_{F, i}(B)$} as
\(\mathcal{B}(B, F \bigsqcup \egraph^i)\).
Equivalently, this is the number of boards
$\sigma \subseteq B$ with $G(\sigma) \cong F \bigsqcup \egraph^i \bigsqcup R$ (where \(R\) is a graph with no edges of the appropriate number of row and column vertices).
\end{definition}
For example, when
$F$ is the empty graph, $r_{F,i}(B)$ is the standard rook number $r_i(B)$. By
convention, \(r_{F,i}(B) = 0\) if \(i < 0\).
We give a formula for the $q$-rook number $M_d(B,q)$ modulo $(q-1)^2$:
\begin{theorem}
\label{thm:q-rookres}
For a board $B$ and non-negative integer $d$, 
\begin{align*}
M_d(B,q) &\equiv (q-1)\bigl(r_{\ZG,d-2}(B) - r_{\SG,d-2}(B) + r_{\WRG,d-1}(B) + r_{\WCG,d-1}(B)\bigr) \\
&+ r_d(B) \pmod{(q-1)^2}.
\end{align*}
\end{theorem}

\begin{proof}

  By Equation~\eqref{eqn:m-sum-connected} and Proposition~\ref{prop:graph-classes}, we have
  \[
    \mathfrak{m}_d(B,q) \equiv r_d(B)(q-1)^d + (q-1)^{d+1} \sum_{(G,\Rank)} r_{G,d-\Rank}(B) \mathcal{O}_\Rank(G,q) \pmod{(q-1)^{d+2}},\]
  where the sum ranges over \(\{(\ZG,2),(\SG,2),(\WCG,1),(\WRG,1)\}\).
  Using Proposition~\ref{prop:og}, and the
  multiplicativity, we get:
  \begin{equation}
  \label{eqn:q-rookmod}
  \begin{split}
  \mathfrak{m}_d(B,q) &\equiv r_d(B)(q-1)^d + (q-1)^{d+1}\Big(r_{\ZG, d-2}(B) + r_{\WRG, d-1}(B)\\ &+ r_{\WCG, d-1}(B) + (q-2) r_{\SG, d-2}(B)\Big) \pmod{(q-1)^{d+2}}.
  \end{split}
  \end{equation}

  Dividing by $(q-1)^d$ gives the desired equation.
\end{proof}

\begin{remark}
\label{rem:q-rook-linear-nonneg}
    Since $r_{\SG, d-2}(B) \leq r_{\ZG, d-2}(B)$, we have $r_{\ZG, d-2}(B) - r_{\SG, d-2)}(B) + r_{\WRG, d-1}(B) + r_{\WCG, d-1}(B) \geq 0$, and the coefficients of $M_d(B,q)\pmod{(q-1)^2}$ are nonnegative.
\end{remark}

\begin{example}
  For example, in the $B = [2]\times [2]$ board $B$ with $d = 2$, we have
  $r_2(B) = 2, r_{\ZG, 0}(B) = 4, r_{\SG, 0} = 3$, and
  $r_{\WRG,1} = r_{\WCG,1} = 0$, so
\begin{align*}
M_2(B, q) &= (q^2 - 1)(q^2-q)/(q-1)^2 \equiv 2 + 3(q-1) \\
& \equiv r_2(B) + (q-1)(r_{\ZG,0} - r_{\SG,0} + r_{\WRG,1} + r_{\WCG,1})\pmod{(q-1)^2}.
\end{align*}

The linear term of $q-1$ also has nonnegative coefficients.
\end{example}

\subsection{Polynomiality of \texorpdfstring{\(\mathfrak{m}_d(B,q)\)}{m\_d(B,q)} modulo \texorpdfstring{$(q-1)^{d+6}$}{(q-1)\^(d+6)}}
\label{sec:q-hit-poly}
Another immediate consequence of Equation~\ref{eqn:m-sum-connected} is the following:
\begin{theorem}
\label{thm:poly-q6}

For any board $B$, $\mathfrak{m}_d(B,q)$ is polynomial modulo $(q-1)^{d+6}$, and
its coefficients modulo \((q-1)^{d+6}\) are integers.
\end{theorem}

\begin{proof}
  Consider Equation~\eqref{eqn:m-sum-connected} modulo \((q-1)^{d+6}\).
  Let \(G_1, \ldots, G_i\) be the connected components of \(G\), and
  \(\Rank_1, \ldots, \Rank_i\) index a nonzero term of the summation.
  Consider the \(k\)th term of the product,
  \[
    \mathcal{O}_{\Rank_k}(G_k, q)(q-1)^{x_k+y_k-\Rank_k-1}.
  \]
  If either
  of \(x_i\) or \(y_i\) is at least \(7\), then \(x_i+y_i-\Rank_i-1\) is at
  least \(6\), as \(\Rank_i\) can be at most \(\min(x_i,y_i)\).
  In this case the entire term of the summation is congruent to 0.
  Therefore, none of the connected components of \(G\) can have \(7\) or more
  of either row or column vertices.
  We know however by Remark~\ref{rem:fano_minimal} that
  \(\mathcal{O}_{\Rank_k}(G_k,q)\) is polynomial if \(G_k\) has at most \(6\)
  row vertices and at most \(6\) column vertices, and so
  \(\mathfrak{m}_d(B,q)\) is congruent to a polynomial.
\end{proof}

\begin{corollary}
\label{cor:poly-hit-q6}
For a board $B$, $M_d(B,x+1)$ and $H_d(B,x+1)$ are both polynomial modulo \(x^6\), and have integer coefficients.
\end{corollary}

\begin{proof}
    By definition,
    \[
    M_d(B,q) = \mathfrak{m}_d(B,q)/(q-1)^d,
    \]
    so \(M_d(B,q)\) is polynomial with integer coefficients modulo \((q-1)^6\) for any \(d\).
    Then Proposition~\ref{prop:q-hit-direct} shows \(H_d(B,q)\) is too.
\end{proof}

\begin{remark}
\label{rem:k=6maximal}
The number $6$ is the maximal $c$ such that $\mathfrak{m}_d(B,q)$ is always
polynomial modulo \((q-1)^{d+c}\). Consider the Fano board $F$ with $d = 7$.
Using Equation~\eqref{eqn:fano-hit-eq}, the coefficient of $(q-1)^{13}$ depends
on the residue class of $q\pmod{2}$. Therefore $\mathfrak{m}_d(F, q)$ is not
polynomial modulo \((q-1)^{d+7}\).
\end{remark}

\section{Hit numbers modulo \texorpdfstring{$(q-1)^2$}{(q-1)\^2}}
\label{sec:higher-hit}

In this section, we consider the $q$-hit number $H_i(B,q)$.
Motivated by \cite[Conjecture 6.7]{lewis_rook_2020}, we conjecture the following.
\begin{conjecture}
\label{conj:q-poly-all-coeffs}
Let \(B\) be a board, and \(P \in \mathbb{Z}[x]\) be a polynomial of degree
$k-1$ such that \(P(x) = H_i(B, x+1) \pmod{x^k}\) for all \(x\) in an unbounded
subset of \(\mathbb{Z}\). Then, \(P\) has non-negative coefficients.
\end{conjecture}

In this section we verify this conjecture for $k=1$ and $k=2$. By \cite[Prop.
3.3]{lewis_rook_2020}, we know that $H_i(B,x+1) \equiv h_i(B) \pmod{x}$, which
is manifestly non-negative. Next, we find an expression for the coefficient of
\(x\) in $H_i(B,x+1)$ modulo \(x^2\). Because $H_i(B,x)$ is a polynomial with
integer coefficients modulo $x^2$ (by Corollary~\ref{cor:poly-hit-q6}), we know
that the coefficient does not depend on \(x\) and is an integer.

Therefore, let \defn{\(C_i(B)\)} be defined by
\begin{equation} \label{eq:def C}
H_i(B, x+1) \equiv C_i(B)x + h_i(B) \pmod{x^2}.
\end{equation}
We seek a formula for \(C_i(B)\) by reducing Equation~\eqref{eqn:q-hitdef}
modulo $(q-1)^2$ and then extracting the coefficient of $t^i$ from both sides.
We show that the coefficient \(C_i(B)\) is non-negative for any board
\(B\).

\subsection{Finding \texorpdfstring{$H_d(B,x+1)$}{H\_d(B,x+1)} modulo \texorpdfstring{$x^2$}{x^2}}
\label{sec:q-hit-mod}

In this subsection, we transfer the results of Section \ref{sec:orbit-q-rook} to
the computation of \(q\)-hit numbers. We use our formula~\eqref{eqn:q-rookmod}
to produce a formula for $H_d(B,x+1) \pmod{x^2}$.

We make the following definitions of generalized hit numbers, analogous to
the generalized rook numbers in Section~\ref{sec:orbit-q-rook}.

\begin{definition}
  For a board $B \subseteq [m]\times[n]$, and a bi-colored graph \(F\) with
  \(x\) row and \(y\) column vertices, we define \defn{$h_{F, d}(B)$} as the
  number of boards $\sigma \subseteq[m]\times[n]$ with
  $G(\sigma) \cong F \bigsqcup \egraph^{\min(m-x, n-y)} \bigsqcup R$, such that
  \(G(\sigma \cap B) \cong F \bigsqcup \egraph^{d} \bigsqcup R'\), where \(R\) and
  \(R'\) are (bi-colored) graphs with no edges.
\end{definition}

Note that, if $F$ is an empty graph, then $h_{F,d}(B)$ is exactly the usual hit number
\(h_d(B)\).

Now we can state the formula which is the main result of this subsection.
\begin{theorem}
\label{thm:q-hitres}
The \(q\)-hit number satisfies
\[
H_d(B,x+1) \equiv C_d(B)x + h_d(B) \pmod{x^2},
\]
where \(C_d(B)\) is given by one of the following formulas.

For boards $B\subseteq [m]\times [n]$ with $m < n$,
\begin{align*}
    C_d(B) &= h_{\ZG, d-2}(B) - h_{\SG, d-2}(B) +  (n-m+1)h_{\WRG, d-1}(B)\\
    &+ \frac{n-d}{n-m} h_{\WCG, d-1}(B) - 2h_{\ZG, d-1}(B) + 2h_{\SG, d-1}(B)
- (n-m+1) h_{\WRG, d}(B)\\
&+ \frac{2d-n-1}{n-m} h_{\WCG, d}(B) + h_{\ZG, d}(B) - h _{\SG, d}(B) + \frac{d-1}{n-m} h_{\WCG, d+1}(B) \\
&+ \frac14\Big(h_d(B)(m-d)(m+2n+d-3) \\
&+ h_{d+1}(B)(2d+2)(n-1) + h_{d+2}(B)(d+2)(d+1)\Big).
\end{align*}
For boards $B \subseteq [n]\times [n]$,
\begin{equation}
\begin{split}
C_d(B) &= h_{\ZG, d-2}(B) - h_{\SG, d-2}(B) + h_{\WRG, d-1}(B) + h_{\WCG, d-1}(B)\\
&- 2h_{\ZG, d-1}(B) + 2h_{\SG, d-1}(B)
- h_{\WRG, d}(B) - h_{\WCG, d}(B)\\
&+ h_{\ZG, d}(B) - h_{\SG, d}(B) + \frac14 \big(h_d(B)(n-d)(3n + d - 3)\\
&+ h_{d+1}(B)(2d+2)(n-1) + h_{d+2}(B) (d+2)(d+1)\big).
\end{split}
\end{equation}
\end{theorem}

Recall that \(C_d(B)\) is an integer by Corollary~\ref{cor:poly-hit-q6}. In Section~\ref{sec:coeff-hit-positive} we will show $C_d(B)$ is \emph{non-negative}.

We start with an example, to make the notation clear.

\begin{example}
Let $B = [2]\times [2]$ be the full \(2 \times 2\) board. By direct computation,
we have

\[H_2(B,x+1) = x^2 + 3x + 2 \equiv 3x + 2\pmod{x^2}.\]

Indeed, $h_{\ZG,0}(B) = \mathcal{B}(B,\ZG) = 4$ and
$h_{\SG, 0}(B) = \mathcal{B}(B,\SG) = 1$, so
$h_{\ZG, 0}(B) - h_{\SG, 0}(B) = C_2(B) = 3$.
\end{example}

Before we begin the proof, we present the key ingredients of the computation. The following results about \(q\)-analogues have standard proofs, which we omit:
\begin{lemma}
\label{lem:q-number-linear-coeffs}
For all integers \(a, n \geq 1\), we have the following congruences of polynomials:
\begin{align*}
q^n &\equiv n(q-1) +1 \pmod{(q-1)^2},\\
(n)_q &\equiv \binom{n}{2}(q-1) + n \pmod{(q-1)^2},\\
[n]!_q &\equiv \frac{n!}{2}\binom{n}{2}(q-1) + n!\pmod{(q-1)^2}\\
\qbinom{n}{i}_q &\equiv \binom{n}{i} \frac{i(n-i)}{2}(q-1)  + \binom{n}{i} \pmod{(q-1)^2}\\
(a;q)_n &\equiv  \left(-\binom{n}{2}(1-a)^{n-1} a\right)(q-1) +  (1-a)^{n} \pmod{(q-1)^2}.
\end{align*}
\end{lemma}

We use the following relation between the classical rook and hit numbers.
\begin{lemma}
\label{lem:rookhit-fallingfactweight}
For a board $B$ and fixed $k$, we have
\begin{align*}
&\sum_{i=0}^m r_i(B) \frac{(n-i)!}{(n-m)!} i(i-1)\ldots (i-k + 1) (t-1)^i \\
&= \sum_{i=0}^{m}i(i-1)\ldots (i-k+1)(t-1)^kt^{i-k}h_i(B).
\end{align*}
\end{lemma}
\begin{proof}
Take \(k\) derivatives of the equation~\eqref{eqn:hit-rook-relation} and then multiply by $(t-1)^k$.
\end{proof}

Finally, we have generalized rook-hit relations, analogous to the classical
rook-hit relation from Equation~\eqref{eqn:hit-rook-relation}.

Let \(B \subseteq [m]\times[n]\) be a board, and \(G\) be a bi-colored graph with \(x\)
row vertices and \(y\) column vertices. Let \(\Min = \min(m-x,n-y)\), and
\(\Max = \max(m-x, n-y)\).

\begin{proposition}
  \label{prop:generalized-hit-rook-relation}
  The generalized rook numbers are related to the generalized hit numbers by the equation
  \[\sum_{i=0}^\Min r_{G,i}(B) \frac{(\Max-i)!}{(\Max-\Min)!} (t-1)^i = \sum_{i=0}^\Min h_{G,i}(B) t^i.\]
\end{proposition}
\begin{proof}
We define a generating function for certain arrangements of boards, and equate two different formulas for it.
\begin{definition}
  A \emph{hit pair} of weight \(i\) is a pair of boards
  \(\omega \subseteq \sigma \subseteq [m]\times [n]\) such that
  \[
    \omega \subseteq B \text{ and } G(\omega) \cong G \bigsqcup \egraph^{i}
  \]
  and
  \[
    G(\sigma) \cong G \bigsqcup \egraph^\Min.
  \]
  We write \(\wt((\omega,\sigma)) = \wt(\omega) = i\).
\end{definition}

Let \(T(t)\) be the (shifted) generating function for hit pairs by weight:
\[
  T(t) = \sum_{(\omega, \sigma)} (t-1)^{\wt(\omega)}.
\]

If we fix \(\omega\), the number of \(\sigma\) forming a hit pair with
\(\omega\) is \(r_{\Min-i}([m-x-i]\times [n-y-i])\) independent of which board
\(\omega\) is. Since the smaller of \(m-x-i\) and \(n-y-i\) is equal to
\(\Min-i\), and the bigger is equal to \(\Max-i\),
\[
  r_{\Min-i}([m-x-i]\times[n-y-i]) = (\Max-i)_{\Min-i} = \frac{(\Max - i)!}{(\Max - \Min)!}.
\]
The number of choices of \(\omega\) with weight \(i\) is exactly the rook number \(r_{G,i}(B)\), and so
\[
  T(t) = \sum_i r_{G,i}(B)\frac{(\Max-i)!}{(\Max-\Min)!}(t-1)^i.
\]

Now, the number of \(\sigma\) such that
\(G(\sigma \cap B) \cong G \bigsqcup \egraph^j\) is the hit number \(h_{G,j}(B)\). If
we fix such a \(\sigma\), then the number of boards \(\omega \subseteq \sigma\)
forming a hit pair of weight \(i\) with \(\sigma\) is \(\binom{j}{i}\). Thus,
\begin{align*}
  T(t) &= \sum_{j=0}^m h_{G,j}(B)\sum_{i=0}^j \binom{j}{i}(t-1)^i \\
  &= \sum_{j=0}^m h_{G,j}(B)t^j.
\end{align*}
\end{proof}

We are now ready to complete the proof of Theorem~\ref{thm:q-hitres}.
\begin{proof}[Proof of Theorem~\ref{thm:q-hitres}]

Consider Equation~\eqref{eqn:q-hitdef}, reproduced here:
\[
  \sum_{i=0}^{m} H_i(B,q) t^i = q^{\binom{m}{2}} \sum_{i=0}^{m} M_i(B,q) \frac{[n-i]!_q}{[n-m]!_q}(-1)^i (t;q^{-1})_i.
\]

We wish to reduce the right-hand side modulo $(q-1)^2$ then extract the
coefficients for $t^i$ modulo $(q-1)^2$.

By Theorem~\ref{thm:q-rookres},
\[
  M_i(B,q) \equiv \left(r_{\ZG,i-2}(B) - r_{\SG,i-2}(B) + r_{\WRG,i-1}(B) + r_{\WCG,i-1}(B)\right)(q-1) + r_i(B) \pmod{(q-1)^2}
\]
By Lemma~\ref{lem:q-number-linear-coeffs},
\[
  \frac{[n-i]!_q}{[n-m]!_q} \equiv
  \frac{(n-i)!}{(n-m)!}\frac{(m-i)(2n-m-i-1)}{4}(q-1) + \frac{(n-i)!}{(n-m)!} \pmod{(q-1)^2},
\]
and
\[
  (-1)^i(t;q^{-1})_i \equiv -\binom{i}{2}t(t-1)^{i-1}(q-1) + (t-1)^i \pmod{(q-1)^2}.
\]
We distribute and sum over all \(i\).
\begin{align*}
  \sum_{i=0}^m H_i(B,q) t^i \equiv \Big[ & (q-1)\sum_{i=0}^m \left(r_{\ZG,i-2}(B) - r_{\SG,i-2}(B) + r_{\WRG,i-1}(B) + r_{\WCG,i-1}(B)\right)\frac{(n-i)!}{(n-m)!}(t-1)^i \\
                                                    +& (q-1)\sum_{i=0}^m r_i(B)\frac{(n-i)!}{(n-m)!}\frac{(m-i)(2n-m-i-1)}{4}(t-1)^i \\
  +& (q-1)\sum_{i=0}^m r_i(B)\frac{(n-i)!}{(n-m)!}\left(-\binom{i}{2}t(t-1)^{i-1}\right) \\
  +& \sum_{i=0}^m r_i(B)\frac{(n-i)!}{(n-m)!}(t-1)^i\Big]q^{\binom{m}{2}} \pmod{(q-1)^2}.
\end{align*}

Next we rewrite in terms of (generalized) hit numbers. The constant coefficient
is
\[
  \sum_{i=0}^m r_i(B)\frac{(n-i)!}{(n-m)!}(t-1)^i = \sum_{i=0}^m h_i(B)t^i,
\]
by Equation~\ref{eqn:hit-rook-relation}.

The coefficient of \((q-1)\) has three terms. The two sums involving classical rook
numbers can be rewritten using Lemma~\ref{lem:rookhit-fallingfactweight}:

\begin{align*}
  & \sum_{i=0}^m r_i(B)\frac{(n-i)!}{(n-m)!}\frac{(m-i)(2n-m-i-1)}{4}(t-1)^i \\
  = \frac{1}{4}& \sum_{i=0}^m h_i(B)t^{i-2}\left((2n-m-1)mt^{2} - (2n-2)i(t-1)t + 2(t-1)^2\binom{i}{2}\right),
\end{align*}
and
\[
  \sum_{i=0}^m r_i(B)\frac{(n-i)!}{(n-m)!}\left(-\binom{i}{2}t(t-1)^{i-1}\right) = - \sum_{i=0}^m h_i(B) \binom{i}{2} t^{i-1}(t-1).
\]

The first sum can be rewritten using
Proposition~\ref{prop:generalized-hit-rook-relation} for each of the four
isomorphism classes of bi-colored graph. Each term in the sum
corresponds to a type of orbit, and looks like
\[
  f_{G,\Rank}(B,t) := \sum_{i=\Rank}^{\Min+\Rank} r_{G,i-\Rank}(B)\frac{(n-i)!}{(n-m)!}(t-1)^i,
\]
where \(\Rank\) is the rank of the submatrix with graph \(G\) for all matrices
in those orbits.

Some difficulties arise when the
quantity \((n-i)!\) does not coincide with the quantity \((\Max + \Rank - i)!\) in the
relation of Proposition~\ref{prop:generalized-hit-rook-relation}, which is only the case for the \(\WCG\) class--for now, we encapsulate the details in the following Lemma.
\begin{lemma}
  \label{lem:hit-rook-qterm}
  The function
\begin{equation*}
f_{G,\Rank}(B,t) = \sum_{i=\Rank}^{\Min+\Rank} r_{G,i-\Rank}(B)\frac{(n-i)!}{(n-m)!}(t-1)^i,
\end{equation*}
may also be written as
\begin{equation}
  \label{eqn:hit-rook-qterm}
  f_{G,\Rank}(B,t) = (t-1)^\Rank \frac{(\Max-\Min)!}{(n-m)!} (-1)^{\Defect} \sum_{i=0}^{\Min} \left[ \sum_{\ell = 0}^{\Defect} h_{G,i+\ell}(B)(-1)^\ell\sum_{j=\ell}^{\ell+i}\binom{\Defect}{j}\binom{j}{\ell} (-\Max-1)_{\Defect-j} (i+\ell)_{j}\right]t^i,
\end{equation}
where \(\Defect := n-\Rank-\Max\) is called the \defn{defect} of the pair \((G,\Rank)\)
on the board \(B\).
\end{lemma}

The proof is by combining defect-many derivatives of
Proposition~\ref{prop:generalized-hit-rook-relation}.
Note that \(\Max, \Min,\) and \(\Defect\) depend on the dimensions of \(G\) and
the rank \(\Rank\). We remark that for \(3\) of the \(4\) \((G,\Rank)\) pairs,
the defect is zero, and
Equation~\eqref{eqn:hit-rook-qterm} simplifies to
\[
  f_{G,\Rank}(B,t) = (t-1)^\Rank \frac{(\Max-\Min)!}{(n-m)!} \sum_{i=0}^\Min h_{G,i}(B)t^i.
\]

We return now to the task of rewriting
the sum
\[
  \sum_{i=0}^m \left(r_{\ZG,i-2}(B) - r_{\SG,i-2}(B) + r_{\WRG,i-1}(B) + r_{\WCG,i-1}(B)\right)\frac{(n-i)!}{(n-m)!}(t-1)^i,
\]
which is \(f_{\ZG,2}(B,t) - f_{\SG,2}(B,t) + f_{\WRG,1}(B,t) + f_{\WCG,1}(B,t)\).
Applying Lemma~\ref{lem:hit-rook-qterm}, we have
\begin{align*}
  & f_{\ZG,2}(B,t) - f_{\SG,2}(B,t) + f_{\WRG,1}(B,t) + f_{\WCG,1}(B,t) \\
  =& (t-1)^2\sum_{i=0}^{m-2} [h_{\ZG,i}(B) - h_{\SG,i}(B)]t^i\\
     +& (t-1)\sum_{i=0}^{m-2} \left[\frac{(n-m-1)!}{(n-m)!}h_{\WRG,i}(B)t^i\right] \\
    +& f_{\WCG,1}(B,t).
\end{align*}
The \(\WCG\) term can be simplified drastically by dividing into two cases, depending on whether
\(m = n\) or \(m < n\). If \(m = n\), then \(\Max = m-1\), \(\Min = m-2\), and
\(\Defect = 0\).
Then
\[
  f_{\WCG,1}(B,t) = (t-1)\sum_{i=0}^{m-2} h_{\WCG,i}(B)t^i.
\]

If \(m < n\), \(\Max = n-2\), \(\Min = m-1\), and \(\Defect = 1\). Then we have
a more complicated formula:
\[
  f_{\WCG,1}(B,t) = -(t-1)\sum_{i=0}^{m-1} \left[
    ih_{\WCG,i}(B)
    -
    (i+1)h_{\WCG,i+1}(B) \right]t^i.
\]

Putting everything together, we have (when \(m < n\))
\begin{align*}
& \sum_{i=0}^m \left(r_{\ZG,i-2}(B) - r_{\SG,i-2}(B) + r_{\WRG,i-1}(B) + r_{\WCG,i-1}(B)\right)\frac{(n-i)!}{(n-m)!}(t-1)^i\\
=& \sum_{i=0}^m t^{i-2}(t-1) \Big[(t-1)\big(h_{\ZG,i-2}(B) - h_{\SG,i-2}(B) - \frac{i-1}{n-m}h_{\WCG,i-1}(B)\big)\\
+& t\big((n-m+1)h_{\WRG,i-1}(B) + \frac{n-1}{n-m}h_{\WCG,i-1}(B)\big)\Big].
\end{align*} 

For $m = n$, we get
\begin{align*}
& \sum_{i=0}^n \left(r_{\ZG,i-2}(B) - r_{\SG,i-2}(B) + r_{\WRG,i-1}(B) + r_{\WCG,i-1}(B)\right)(n-i)!(t-1)^i\\
=& \sum_{i=0}^n t^{i-2}(t-1)\Big[ (t-1)(h_{\ZG, i-2}(B) - h_{\SG, i-2}(B)) + t(h_{\WRG,i-1} + h_{\WCG,i-1})\Big]
\end{align*}

We extract the coefficient of $t^i$:

When \(m < n\), we obtain
\begin{align*}
    C_i(B) &= h_{\ZG, i-2}(B) - h_{\SG, i-2}(B) +  (n-m+1)h_{\WRG, i-1}(B) \\
    &+ \frac{n-i}{n-m} h_{\WCG, i-1}(B)- 2h_{\ZG, i-1}(B) + 2h_{\SG, i-1}(B) \\
&- (n-m+1) h_{\WRG, i}(B) + \frac{2i-n-1}{n-m} h_{\WCG, i}(B)\\
&+ h_{\ZG, i}(B) - h _{\SG, i}(B) + \frac{i-1}{n-m} h_{\WCG, i+1}(B) \\
&+ \frac14\Big(h_i(B)(m-i)(m+2n+i-3) \\
&+ h_{i+1}(B)(2i+2)(n-1) + h_{i+2}(B)(i+2)(i+1)\Big).
\end{align*}

When \(m = n\), it simplifies to
\begin{align*}
C_i(B) &= h_{\ZG, i-2}(B) - h_{\SG, i-2}(B) + h_{\WRG, i-1}(B) + h_{\WCG, i-1}(B)\\
&- 2h_{\ZG, i-1}(B) + 2h_{\SG, i-1}(B)- h_{\WRG, i}(B) - h_{\WCG, i}(B) \\
&+ h_{\ZG, i}(B) - h_{\SG, i}(B) + \frac14 \big(h_i(B)(n-i)(3n + i - 3) \\
&+ h_{i+1}(B)(2i+2)(n-1) + h_{i+2}(B) (i+2)(i+1)\big).
\end{align*}
\end{proof}

\subsection{Positivity in q-hit number coefficients}
\label{sec:coeff-hit-positive}
In Theorem~\ref{thm:q-hitres} we gave an alternating formula for \(C_i(B)\) defined in \eqref{eq:def C}. In this section, we show \(C_i(B) \geq 0\). This verifies Conjecture~\ref{conj:q-poly-all-coeffs} for \(k=2\), as $H_i(B,q) \equiv (q-1)C_i(B) + h_i(B) \pmod{(q-1)^2}$, which has nonnegative coefficients.

\begin{example}
  \label{exa:derangements}
Let \(B=\{(1,1),\ldots,(n,n)\}\) be the support of the \(n \times n\) identity matrix. The \(q\)-rook
numbers of its complement \(\overline{B}\) have been studied (\cite{lewis_matrices_2011},\cite{ravagnani_codes_2015}) as a \(q\)-analogue of derangements.

The \(q\)-hit numbers of any board and its complement satisfy a reciprocity relation
\cite[Prop. 3.9]{lewis_rook_2020}
\[
  H_d(\overline{B},q) = q^{n(n-d)-\#B}H_{n-d}(B,q),
\]
so we have, from \cite[Prop 3.5]{lewis_rook_2020},
\[H_{n-d}(\overline{B}, q) = q^{n(d-1)}q^{\binom{d+1}{2} + \binom{n}{2}} \sum_{j=d}^{n} M_j(B,q) [n-j]!_q \qbinom{j}{d}_q (-1)^{j+d} q^{-j(d)}\]
\[= q^{n(d-1) + \binom{d+1}{2} + \binom{n}{2}} \sum_{k=0}^{n-d} \binom{n}{k+d} [n-k-d]!_q [k+d]!_q [d]!_q^{-1} [k]!_q^{-1} (-1)^{k} q^{-(k+d)d}\]
\[= q^{n(d-1) - \binom{d}{2} + \binom{n}{2}} \sum_{k=0}^{n-d} (-1)^{k} \binom{n}{k+d} [n-d-k!]_q \qbinom{d+k}{d}_q.\]
From \ref{lem:q-number-linear-coeffs}, we get
\begin{align*}
H_{n-d} &\equiv \left((q-1)\left(\binom{n}{2} - \binom{d}{2} + n(d-1)\right) + 1\right) \sum_{k=0}^{n-d} \Bigg[(-1)^{k}\binom{n}{k+d} (n-d-k)!\\
&\cdot \left(\frac12 \binom{n-d-k}{2}(q-1) + 1\right) \binom{d+k}{d} \left(\frac{dk}{2}(q-1) + 1\right)\Bigg] \pmod{(q-1)^2}.
\end{align*}
Therefore, we have
\[C_{n-d}(\overline{B}) = \sum_{k=0}^{n-d} \frac{n!}{k!d!} (-1)^{k} \left( \binom{n}{2} - \binom{d}{2} + n(d-1) + \frac12 \binom{n-d-k}{2} + \frac{dk}{2}\right).\]
This is nonnegative, since for $k = 2j$, the term inside the summation is greater in magnitude than for $k = 2j+1$, and the inside term is positive for $k = 2j$ and negative for $k = 2j+1$.

\end{example}

\begin{theorem}
\label{thm:q-coeff1-pos}
For a board $B \subseteq [n]\times[n]$, the linear $q$-hit coefficient $C_i(B)$  is a non-negative integer.
\end{theorem}

We prove that \(C_i(B)\) is positive using a series of inequalities relating our generalized hit numbers to the usual ones.

\begin{lemma}
For boards $B\subseteq[n]\times[n]$, we have
\begin{equation}
\label{ineqn:square-chain-1}
 h_{\WRG,i-1}(B) + h_{\WCG,i-1}(B) - 2h_{\ZG,i-1}(B) + 2h_{\SG,i-1}(B) + \frac14 (2i+2)i h_{i+1}(B) \geq 0
\end{equation}
and
\begin{equation}
\label{ineqn:square-chain-2}
\begin{split}
&h_{\ZG,i}(B) - h_{\SG,i}(B) - h_{\WRG,i}(B) - h_{\WCG,i}(B)+\\ &\frac14 (2i+2)(n-i-1) h_{i+1}(B) + \frac14 (i+1)(i+2)h_{i+2}(B) \geq 0.    
\end{split}
\end{equation}
\end{lemma}
\begin{proof}

Define a \defn{square-chain} as a set of cells \(t \subset [n] \times [n]\) whose associated graph \(G(t)\) consists of the union of $n-2$ disjoint edges and a $K_{2,2}$. For a square-chain $t$, let $T(t)$ be the set of cells in \(t\) corresponding to the \(K_{2,2}\). For a board $B \subseteq [n]\times[n]$, define $S_i(B)$ as the set of square chains $t$ such that $\#((t\backslash T(t)) \cap B) = i - 2$. 

For each board $\omega\subseteq[n]\times[n]$ with $G(\omega) \cong \ZG\bigsqcup \egraph^{n-2}$ and $G(\omega \cap B) \cong \ZG\bigsqcup \egraph^{i-2}$, there is exactly one $t \in S_i(B)$ such that $\omega \subseteq t$, and there are no $t \in S_k(B)$ for $i\neq k$ such that $\omega \subseteq t$. The same is true mutatis mutandis for the isomorphism classes \(\SG, \WRG,\) and \(\WCG\).

For board $\omega \subseteq[n]\times[n]$ with $G(\omega)$ consisting of $n$ disjoint edges and $G(\omega\cap B)$ consisting of $i$ disjoint edges, there are exactly $\binom{i}{2}$ different $t \in S_i(B)$ such that $\omega \subseteq t$. This is because we choose two cells $c_1, c_2 \in \omega \cap B$ that are part of $T(t)$ in $\binom{i}{2}$ ways, and such a choice fixes the other two cells in \(T(t)\). Similarly, there are exactly $i(n-i)$ such $t \in S_{i+1}(B)$ such that $\omega \subseteq t$.

For a square-chain \(t\), let $C_{t, \ZG, i}(B)$ be defined as the number of boards $\omega \subseteq[n]\times[n]$ with $G(\omega) \cong \ZG\bigsqcup \egraph^{n-2}$ such that $G(\omega \cap B) \cong \ZG\bigsqcup \egraph^{i}$ and $\omega \subseteq t$. Define $C_{t, \SG,i}(B), C_{t, \WRG,i}(B), C_{t, \WCG,i}(B)$ similarly. Also, define $C_{t, i}(B)$ as the number of boards $\omega\subseteq [n]\times[n]$, where $G(\omega)$ consists of $n$ disjoint edges such that $\omega \subseteq t$ and $\#(\omega\cap B) = i$.

A choice of \(\sigma \subset [n] \times [n]\) such that
\(G(\sigma) \cong \ZG\bigsqcup \egraph^{i-2}\) fixes a unique choice of \(t\) in
\(S_i(B)\), and so
\[\sum_{t\in S_i(B)} C_{t,\ZG,i-2}(B) = h_{\ZG,i-2}(B).\]
The same is true mutatis mutandis for $\SG, \WRG$, and $\WCG$.

Similarly, for the usual hit numbers, we have the relations
\[\sum_{t\in S_i(B)} C_{t, i}(B) = \binom{i}{2} h_i(B)\]
and
\[\sum_{t\in S_{i+1}(B)} C_{t,i}(B) = i(n-i) h_i(B).\]
Observe that if $t\in S_i(B)$, then the values $C_{t,\ZG,i-2}(B), C_{t,\SG,i-2}(B), C_{t,\WRG,i-2}(B),$ $C_{t,\WCG,i-2}(B),$ $ C_{t,i}(B),$ $C_{t,i-1}(B),$ $C_{t,i-2}(B)$ can be determined by the intersection \(T(t) \cap B\). 

First, for $t\in S_{i+1}(B)$, we have the inequality
\[C_{t,\WRG,i-1}(B) + C_{t,\WCG,i-1}(B) - 2C_{t,\ZG,i-1}(B) + 2C_{t,\SG,i-1}(B) + C_{t,i+1}(B) \geq0.\]
Since the values can be determined by $T(t)\cap B$, and \(T(t)\) only has \(2^4\) subsets, we can check that the inequality holds for each of the $2^4$ configurations of $T(t)\cap B$, and conclude that it holds for all $t$. Now, we have
\[\sum_{t\in S_{i+1}(B)} C_{t,\WRG,i-1}(B) + C_{t,\WCG,i-1}(B) - 2C_{t,\ZG,i-1}(B) + 2C_{t,\SG,i-1}(B) + C_{t,i+1}(B) \geq0,\]
which is equivalent to
\[ h_{\WRG,i-1}(B) + h_{\WCG,i-1}(B) - 2h_{\ZG,i-1}(B) + 2h_{\SG,i-1}(B) + \frac14 (2i+2)i h_{i+1}(B) \geq 0.\]

Next, for $t \in S_{i+2}(B)$, we have the inequality
\[C_{t,\ZG,i}(B) - C_{t,\SG,i}(B) - C_{t,\WRG,i}(B) - C_{t, \WCG,i}(B) + \frac12 C_{t,i+1}(B) + \frac12 C_{t,i+2}(B) \geq 0.\]
Again, the inequality holds for each of the $2^4$ configurations of $T(t)\cap B$, so this inequality holds for all $t$. Now, we have
\[\sum_{t \in S_{i+2}(B)}C_{t,\ZG,i}(B) - C_{t,\SG,i}(B) - C_{t,\WRG,i}(B) - C_{t, \WCG,i}(B) + \frac12 C_{t,i+1}(B) + \frac12 C_{t,i+2}(B) \geq 0,\]
which is equivalent to
\begin{equation*}
\begin{split}
h_{\ZG,i}(B) - h_{\SG,i}(B) - h_{\WRG,i}(B) - h_{\WCG,i}(B)\\ + \frac14 (2i+2)(n-i-1) h_{i+1}(B) + \frac14 (i+1)(i+2)h_{i+2}(B) \geq 0.    
\end{split}
\end{equation*}
\end{proof}
Now we are ready to complete the proof of Theorem~\ref{thm:q-coeff1-pos}.
\begin{proof}[Proof of Theorem~\ref{thm:q-coeff1-pos}]
Adding inequality~\eqref{ineqn:square-chain-1} and
inequality~\eqref{ineqn:square-chain-2} to the two inequalities
$\frac14 (n-i)(3n+i-3) h_{i}(B)\geq 0$ and
$h_{\ZG,i-2}(B) - h_{\SG,i-2}(B) \geq 0$, we get
\begin{align*}
C_i(B) = &h_{\ZG,i-2}(B) - h_{\SG,i-2}(B) + \frac14 (n-i)(3n+i-3)h_i(B)\\
&+ h_{\WRG,i-1}(B) + h_{\WCG,i-1}(B)
 - 2h_{\ZG,i-1}(B) + 2h_{\SG,i-1}(B)\\
 &+ h_{\ZG,i}(B) - h_{\SG,i}(B) -  h_{\WRG,i}(B) - h_{\WCG,i}(B)\\
+ &\frac14 (2i+2)(n-i-1)h_{i+1}(B) + \frac14 (i+1)(i+2)h_{i+2}(B) \geq 0.
\end{align*}
\end{proof}

\section{Final Remarks}

Since the \(q\)-rook and \(q\)-hit numbers are not necessarily polynomial in
\(q\), it is curious that they seem to have non-negative coefficients (of
\(x = q-1\)) in the sense of Conjecture~\ref{conj:q-poly-all-coeffs}.

We are able to obtain alternating formulas for the coefficient of \(x\) through
the orbit-counting arguments of Section~\ref{sec:orbits}, and in principle this
approach extends further.

Dividing the sum in Equation~\eqref{eqn:m-sum-connected} into parts
corresponding to the finite set of bi-colored graphs and ranks which
can appear, and
analyzing the orbit counts \(\mathcal{O}_{\Rank_k}(G_k,q)\) will yield an
alternating formula even for higher coefficients, but there are two obstructions.
The first is
that the number of terms becomes extremely large, even for the quadratic
coefficient. Particularly for \(q\)-hit numbers, our expansion in
Lemma~\ref{lem:hit-rook-qterm} simplified considerably because the defect can
only be \(0\) or \(1\). The corresponding term in the quadratic coefficient
of \(H_d(B,x+1)\) has many more cases.

The second, more important obstruction is that
\(\mathcal{O}_\Rank(G,q)\) may become non-polynomial as soon as \(G\) contains
at least \(7\) row and \(7\) column vertices, and so one will need to analyze
all polynomials agreeing with the orbit count on unbounded subsets of
\(\mathbb{Z}\). Understanding these orbit counts is expected to be very hard for
general boards.

Conversely, though, by Theorem~\ref{thm:poly-q6}, we know $M_d(B,q)$ modulo $(q-1)^6$ is
an integer polynomial in $q$, and therefore there must be a polynomial formula for
\(M_d(B,x+1) \pmod{x^6}\) in terms of generalized rook numbers
\(\mathcal{B}(B,G)\) for a finite set of graphs \(G\), even if it is too large to write down.
The same is true of \(H_d(B,x+1)\) modulo \(x^6\).

Regarding positivity, the inequalities~\ref{ineqn:square-chain-1} and \ref{ineqn:square-chain-2} have
analogues for generalized hit numbers of other graphs, and we
are curious to see if those analogues could show non-negativity of higher
coefficients than in Theorem~\ref{thm:q-coeff1-pos}.

Of course, showing the coefficients are non-negative leads us to wonder if there
is some nice combinatorial object that the coefficients count.
In Example~\ref{exa:derangements} we see an infinite family of
cases where the \(q\)-hit number and its first coefficient both have nice
counting interpretations.
For instance, \(H_0(\overline{B},q)\) counts arbitrary matrices on the board
\(\overline{B}\), and \(C_0(\overline{B})\) counts cells of the board
\(\overline{B}\). We have so far been unable to find this interpretation in our alternating formula.

\section*{Acknowledgements}

We would like to sincerely thank Alejandro Morales for proposing to study
\cite[Conj. 6.7]{lewis_rook_2020} and thank Joel Lewis, Alejandro
Morales, Tanya Khovanova, and Felix Gotti for many helpful comments and advice.
The authors are both extremely grateful to the MIT PRIMES program for the
opportunity to conduct research together and for their invaluable assistance. JS
was partially supported by NSF grants DMS-1855536 and DMS-22030407. This work was written while JC was a student at the University of Chicago Laboratory Schools.

\printbibliography[
heading=bibintoc,
title={References}
]

\end{document}